\documentclass[a4paper,10pt]{article}
\usepackage{mathtext}
\usepackage[T2A]{fontenc}
\usepackage[utf8]{inputenc}
\usepackage[english]{babel}
\usepackage{latexsym,amsfonts,amssymb,amsmath,longtable,amsthm}
\usepackage{cite}
\usepackage[pic]{xy}
\usepackage{bbm,makeidx}
\usepackage[centerlast,small]{caption}
\usepackage{amsmath}
\title{Twisted conjugacy classes in special and general linear groups}
\author{T.~R.Nasybullov}

\begin{document}

\newcounter{thelem}
\newtheorem{lem}{{\scshape Lemma}}
\newtheorem{ttt}{{\scshape Theorem}}

\maketitle

\begin{center}
\parbox{13cm}{ {\small Abstract.
In this paper we study twisted conjugacy classes and the $R_{\infty}$-property for classical linear groups. In particular, we prove that the general linear group ${\rm GL}_n(K)$ and the special linear group ${\rm SL}_n(K)$ possess $R_{\infty}$-property, if $n\geq3$  and $K$ is an infinite integral domain with trivial group of automorphisms (Theorem 1), or $K$ is an integral domain, which has a zero characteristic and for which ${\rm Aut}(K)$ is torsion (Theorem 2).}}
\end{center}
\section{Introduction}
\par We would like to remind that two elements $x$ and $y$ of a group $G$ are \emph{conjugated}, if there is an element $c$ in $G$, such as $x=cyc^{-1}$. The relation of conjugation is equivalence relation that devides $G$ into conjugacy classes. Relation of twisted conjugation is a generalization of  conjugacy relation, and we would like to remind it. Let $\varphi : G \longrightarrow G$ be an arbitrary automorphism. The elements $x$ and $y$ in $G$ are said to be $ \emph{twisted} ~\varphi$-\emph{conjugated } or simply  $\varphi$-\emph{conjugated} (denoted by $x\sim_{\varphi}y$), if there is an element $c$ in $G$  which satisfies $x=cy\varphi(c^{-1})$. If $\varphi$ is an identical automorphism, then it is simply conjugacy relation. As well as relation of conjugation, the relation of $\varphi$ -- conjugation is equivalence relation and here we can speak about $\varphi$-conjugacy classes. The number $R(\varphi)$ of  $\varphi$-conjugacy classes is called \emph{the Reidemeister number}  of the automorphism $\varphi$. If the Reidemeister number $R(\varphi)$ is infinite for any automorphism $\varphi$ of group $G$, then $G$ is said to have the $R_{\infty}$-property.

The question about the groups that possess $R_{\infty}$-property was formulated by A.~Fel'shtyn and R.~Hill \cite{FH}. This question and other questions related to the twisted conjugacy classes have attracted the attention of many researches \cite{F,FG,FG1,FH,FIT,GW,GW2,KR,BFG}. In particular, A.~Fel'shtyn and D.~Goncalves \cite{FG1} proved that the symplectic group ${\rm Sp}_{2n}(\mathbb{Z})$ possesses $R_{\infty}$-property. The following question arises in this connection: what other linear groups possess $R_{\infty}$-property?

In this paper we prove that the general linear group ${\rm GL}_n(K)$ and the special linear group ${\rm SL}_n(K)$ possess $R_{\infty}$-property, if $n\geq3$  and $K$ is an infinite integral domain with trivial group of automorphisms (Theorem 1), or $K$ is an integral domain, which has a zero characteristic and for which ${\rm Aut}(K)$ is torsion (Theorem 2). The integral domain is understood as a commutative ring with the unit element and without zero divisors.

When the paper was ready to be published, the author got a preprint \cite{MS} from professor P.~Sankaran, where the groups ${\rm SL}_n(\mathbb{Z}), {\rm GL}_n(\mathbb{Z}), {\rm PSL}_n(\mathbb{Z}), {\rm PGL}_n(\mathbb{Z})$ are established to possess $R_{\infty}$-property when $n\geq2$.
\section{Preliminaries}
\par O. O'Meara \cite{M} described all automorphisms of the groups ${\rm GL}_n(K)$ and ${\rm SL}_n(K)$ in the case where $K$ is an integral domain. Let $G$ be one of these groups, then $G$ always possesses the following automorphisms:
\par 1) $\varphi_D$ -- is an automorphism, induced by conjugation  of the any matrix $D$ of ${\rm GL}_n(K)$: \[\varphi_D : A\mapsto D A D^{-1};\]
 If $G={\rm GL}_n(K)$, then it is inner automorphism. We also call it inner automorphism for $G={\rm SL}_n(K)$. All inner automorphisms are known to form a normal subgroup in group of all automorphism.
\par 2) $\bar{\delta}$ -- is a ring automorphism:
\[\bar{\delta} : A=(a_{ij})\mapsto (\delta(a_{ij})),\]  where $\delta$ is an automorphism of the ring $K$;
\par 3) $\Lambda$ -- is a contragradient automorphism:
\[\Lambda : A\mapsto (A^T)^{-1},\] where $^{T}$ means transposition. $\Lambda^2$ is obvious to be an identical automorphism;
\par 4) $\Gamma$ -- is a central automorphism (or a central homothety):
\[\Gamma : A\mapsto \gamma(A) A,\]  where $\gamma$ is a homomorphism of $G$ into its center $Z(G)$. It is well known that $Z(G)$ is the set of scalar matrices. As well as all inner automorphisms, all central automorphisms form a normal subgroup in the group of all automorphisms.

Automorphism, which is a product of automorphisms 1) -- 4),  is called a standart automorphism. In paper \cite{M} it was proved, that any automorphism of groups ${\rm GL}_n(K)$ and ${\rm SL}_n(K)$ is standart if $K$ is an integral domain and $n\geq3$.

We shall prove the auxiliary result to use later.
\begin{lem} Let $K$ be an integral domain, and $M$ be an infinite subset of $K$. Then the set $P=\{f(a): a \in M\}$, where $f$ is an arbitrary polynomial of the power $r\geq 1$, is infinite.
\end{lem}
\textbf{Proof.} If we assume that $P$ is finite, i.e. $P=\{b_1 \dots b_m\}$, then for any element $x$ in $M$ there is a number $i$, such as  $f(x)=b_i$. But this equation can not have more than $r$ roots in an integral domain. Hence $M$ can not have more than $mr$ elements. It contradicts the infinity of $M$.
\begin{flushright}
$\square$
\end{flushright}

We also introduce the following definition to use later.

\textbf{Definition.} \emph{Antitrace} of the matrix $A=\begin{pmatrix}
  a & b \\
  c & d\\
  \end {pmatrix} \in {\rm GL}_2(K)$ is the following element of the \\ring $K$ \[{\rm atr}(A)= b-c. \]

We shall find the property, which is satisfied by an antitrace.
\begin{lem}
Let matrices $X$ and $A$ be in ${\rm GL}_2(K)$. Then the following equality is valid
$${\rm atr}(XAX^T)={\rm atr}(A){\rm det}(X).$$
\end{lem}
\textbf{Proof.}
We note that $A-A^T={\rm atr}(A)\cdot \begin{pmatrix}
   ~~0 & 1 \\
  -1 & 0\\
\end{pmatrix} $, then
\[{\rm atr}(XAX^T)\begin{pmatrix}
   ~~0 & 1 \\
  -1 & 0\\
\end{pmatrix}=XAX^T-XA^TX^T=X(A-A^T)X^T=\] \[={\rm atr}(A)\begin{pmatrix}
   x_{11} & x_{12} \\
  x_{21} & x_{22}\\
\end{pmatrix}\begin{pmatrix}
   ~~0 & 1 \\
  -1 & 0\\
\end{pmatrix}\begin{pmatrix}
   x_{11} & x_{21} \\
  x_{12} & x_{22}\\
\end{pmatrix}={\rm atr}(A)\begin{pmatrix}
   0 & x_{11}x_{22}-x_{12}x_{21} \\
  x_{12}x_{21}-x_{11}x_{22} & 0\\
\end{pmatrix}=\] \[={\rm atr}(A)det(X)\begin{pmatrix}
   ~~0 & 1 \\
  -1 & 0\\
\end{pmatrix}.\]
Hence the result follows.
\begin{flushright}
$\square$
\end{flushright}

\section{ $R_{\infty}$-property in the groups $GL_n(K)$ and $SL_n(K)$ }

One of the main results of this paper is
\begin{ttt} Let $G$ be a general linear group ${\rm GL}_n(K)$ or a special linear group ${\rm SL}_n(K)$ with $n\geq 3$, and $K$ be an infinite integral domain. Let also $\Phi$ be a subgroup of $ {\rm Aut}~G$, which is generated by all inner, central and contragradient automorphisms. Then for any automorphism $\varphi$ in $\Phi$ the number of $\varphi$-conjugacy classes in $G$ is infinite. Particularly, if $K$ is such a ring that $\Phi$ and ${\rm Aut}~G$ are equal, then $G$ possesses $R_{\infty}$-property.
\end{ttt}

\textbf{Proof.}
 Let $\varphi \in \Phi$, then $\varphi=f_1 \dots f_m$, where $f_i$ is either inner, or central or contragradient automorphism. Now using the normality of the subgroup of all inner automorphisms and the subgroup of all central automorphisms in the group of all automorphisms we can rewrite $\varphi$ as:
 \[\varphi = \varphi_D\Lambda^r\Gamma.\]
 Here $\varphi_D$ is inner automorphism, $\Lambda$ is contragradient automorphism and $\Gamma$ is central automorphism of the group $G$.
Now we examine two cases depending on the evenness of $r$.
\par {\it Case} 1. $r$ is an even number. Then $\varphi=\varphi_D \Gamma$. We shall find the condition, which is necessary for $\varphi$-conjugation of matrices $X$ and $Y$ in $G$. By definition, it means that there exists a matrix $Z$ in $G$ for which the following equality is valid \[ X=ZY\varphi(Z^{-1})=ZY\varphi_D\Gamma(Z^{-1})=ZY\varphi_D(\gamma(Z^{-1})Z^{-1})=ZYD\gamma(Z^{-1})Z^{-1}D^{-1}.\]

If we multiply this equality on the right by $D$ then we have the equality:
\begin{equation}
XD=ZYD\gamma(Z^{-1})Z^{-1}.
\end{equation}

To prove that the number of $\varphi$-conjugacy classes is infinite, i.e. $R(\varphi)=\infty$, we need to specify the infinite number of pairwise different $\varphi$-conjugacy classes. According to the established equality (1) it is sufficient to find such a set of matrices $A_1, A_2,...$ in ${\rm SL}_n(K)$, where for any different $i$ and $j$ there is no matrix $C$ in $G$, such as the following equality holds
\[A_iD=CA_jD\gamma(C^{-1})C^{-1}.\]
 Since ${\rm det}~A_i=1$, then if we denote $A_iD=B_i$, we approach the problem of constructing an infinite set of matrices $B_1, B_2, ...$ with given determinants ${\rm det}(B_i)={\rm det}(D)=d$, where for any different $i$ and $j$ there is no matrix $C$ in $G$ which satisfies the condition \[B_i=CB_j\gamma(C^{-1})C^{-1}.\]
  We shall find the condition, when two matrices $X$ and $Y$ with equal non-zero determinants are related by
\begin{equation}
X=ZY\gamma(Z^{-1})Z^{-1}
\end{equation}
for some matrix $Z$ in $G$. If we calculate the determinant in equality (2), we have
\[{\rm det}(X)={\rm det}(\gamma(Z^{-1})){\rm det}(Y)=z^n{\rm det}(Y),\]
where $\gamma(Z^{-1})={\rm diag}(z, ...,z)$ is a scalar matrix.

Hence, using the equality ${\rm det}(X)={\rm det}(Y)$, we have $z^n=1$. But matrices $X$ and $\gamma(Z^{-1})Y$ are conjugated, therefore their traces must be equal, i.e.
\[{\rm tr}(X)={\rm tr}(\gamma(Z^{-1})Y)=z{\rm tr}(Y).\]
 If we raise this equality to the power $n$ using the equality $z^n=1$, we have that the equality $({\rm tr}(X))^n=({\rm tr}(Y))^n$ is necessary for satisfaction of equality (2).

As shown in Lemma 1 the set $P=\{b^n: b \in K\}$ is infinite. Let $\{b_i\}$ be an infinite set of elements of the ring $K$, such that  $b_i^n\neq b_j^n$ when $i\neq j$. It is easy to construct the set of matrices $B_1, B_2, ...$ with given determinant $d$ and with ${\rm tr}(B_i)=b_i$. These are for example the following matrices:
$$
B_i=\begin{pmatrix}
   b_i+2-n & d & 0 \\
   -1 & 0&0\\
   0&0&I_{n-2}
\end{pmatrix}
,$$
  where $I_{n-2}$ is the unitary matrix of $(n-2)\times(n-2)$ size. Then, using the necessary condition, we have that the following equality
\[B_i= CB_j\gamma(C^{-1})C^{-1} \]
does not hold for any matrix $C$ in $G$ when $i\neq j$, i.e. $B_1, B_2, ...$ is required set of matrices. It means, that $R(\varphi)=\infty$ .

\par {\it Case} 2. $r$ is an odd number. Then $\varphi=\varphi_D \Lambda \Gamma$.
 Two matrices $X$ and $Y$ in $G$ are  $\varphi$-conjugated if there exists a matrix $Z$ in $G$, which satisfies the equality
  \begin{multline*}X=ZY\varphi(Z^{-1})=ZY\varphi_D\Lambda \Gamma(Z^{-1})=ZY\varphi_D\Lambda(\gamma(Z^{-1})Z^{-1})=\\
=ZY\varphi_D(\gamma(Z)Z^T)=ZYD\gamma(Z)Z^TD^{-1}.
\end{multline*}
If we multiply it on the right by $D$, then we have
$$
XD=ZYD\gamma(Z)Z^T.
$$
Hence to show that $R(\varphi)=\infty$, it is sufficient to find an infinite set of matrices $A_1, A_2,...$ in ${\rm SL}_n(K)$ where for any different $i$ and $j$ there is no matrix $C$ in $G$ for which the following equality holds
 \[A_iD=CA_jD\gamma(C)C^T.\]
  Then again after denoting $B_i=A_iD$, we have a problem of constructing matrices $B_1, B_2, ...$ with equal determinants, such as for different $i$ and $j$ there is no matrix $C$ in $G$, which satisfies the following equality \begin{equation}B_i=CB_j\gamma(C)C^T.\end{equation}

 We examine the matrices $X$ and $Y$ in ${\rm GL}_2(K)$ with equal determinants and find the necessary condition, when for some matrices $Z$ and $T={\rm diag}(t,t)$ in $G$ the following equality holds
\begin{equation}
X=ZYTZ^T.
\end{equation}
 If we calculate the determinant in equality (4) we have $${\rm det}(X)=({\rm det(Z)})^2{\rm det}(Y){\rm det}(T)=(t\cdot {\rm det}(Z))^2{\rm det}(Y).$$
 In consideration of the equality ${\rm det}(X)={\rm det}(Y)$, we conclude that $(t\cdot {\rm det}(Z))^2=1$. Now we shall calculate an antitrace in equality (4): $${\rm atr}(X)={\rm atr}(ZYTZ^T)={\rm atr}(YT){\rm det}(Z)=t \cdot {\rm atr}(Y){\rm det}(Z).$$
Squaring this equality using the condition $(t \cdot {\rm det}(Z))^2=1$, we conclude that the condition $({\rm atr}(X))^2=({\rm atr}(Y))^2$ is necessary for the equality (4) to be held.

As shown in Lemma 1, the set $P=\{b^2:~b\in K\}$ is infinite. Let $\{b_i\}_{i=1}^{\infty}$ be such an infinite set of elements of the ring $K$, where  $b_i^2\neq b_j^2$ if $i\neq j$. The set
$$
\widehat{B}_i=\begin{pmatrix}
   d & b_i \\
   0 & 1\\
   \end{pmatrix}
$$
is an example of an infinite set of matrices $\widehat{B}_1, \widehat{B}_1, ...$ in $ {\rm GL}_2(K)$ with equal determinants,  such as ${\rm atr}(\widehat{B}_i)=b_i$. Then, according to the necessary condition, when $i \neq j$ the equality
\[\widehat{B}_i= C\widehat{B}_j FC^T\]
can not be carried out for any matrices $C$ and $F={\rm diag}(f,f)$ in ${\rm GL}_2(K)$. Since $({\rm atr}(\widehat{B}_i))^2 \neq ({\rm atr}(\widehat{B}_j))^2$ when $i\neq j$, then it can not be more then one zero element in the set ${\rm atr}(\widehat{B}_i)$. After removing the matrix with zero antitrace, we have that all $\widehat{B}_i$ are nonsymmetric matrices.

Now let ${B}_i=\begin{pmatrix}
      I_{n-2}& 0 \\
     0 & \widehat{B}_i\\
\end{pmatrix}.$ Then $B_1, B_2, \dots$ are stated to be the required set of matrices. Really, supposing matrices $B_i$ and $B_j$ are related by (3) for certain $i\neq j$. It means that there exists matrix  $C=\begin{pmatrix}
   \alpha & \beta \\
   \gamma & \delta \\
\end{pmatrix} \in {\rm GL_n}(K)$, that satisfies the following equality
\[B_i=C\gamma(C)B_jC^T,\]
where $\alpha \in {\rm M}_{(n-2)\times (n-2)}(K)$, $\beta \in {\rm M}_{(n-2) \times 2}(K)$, $\gamma \in {\rm M}_{2\times (n-2)}(K)$, $\delta \in {\rm M}_{2 \times 2}(K)$. We give now a detailed description of the right side of this equality:
\[C\gamma(C)B_jC^T=\begin{pmatrix}
   \alpha & \beta \\
   \gamma & \delta \\
\end{pmatrix}\cdot\gamma(C)\cdot B_j \cdot\begin{pmatrix}
   \alpha^T & \gamma^T \\
   \beta^T & \delta^T \\
\end{pmatrix}=\gamma(C)\cdot\begin{pmatrix}
   \alpha \alpha^T +\beta \widehat{B}_j \beta^T & \alpha \gamma^T +\beta \widehat{B}_j \delta^T \\
   \gamma \alpha^T +\delta \widehat{B}_j \beta_T & \gamma \gamma^T +\delta \widehat{B}_j \delta^T\\
\end{pmatrix}.\]
Let $\gamma(C)= c_1 \cdot I_n$, then the equality $B_i=C\gamma(C)B_jC^T$ implicates a system of equalities:
\begin{gather}
 c_1 I_{n-2}(\alpha \alpha^T +\beta \widehat{B}_j \beta^T)=I_{n-2} ;\\
  \alpha \gamma^T +\beta \widehat{B}_j \delta^T=O_{n-2,2}; \\
 \gamma \alpha^T +\delta \widehat{B}_j \beta_T=O_{2,n-2}; \\
 c_1 I_2(\gamma \gamma^T +\delta \widehat{B}_j \delta^T)=\widehat{B}_i,
\end{gather}
where $O_{k,l}$ is the zero matrix of the $k\times l$ size. If we calculate the antitrace in the equality (8) in consideration of $\gamma \gamma^T$ is a symmetrical matrix, we have
\[{\rm atr}(c_1 I_2(\gamma \gamma^T+\delta \widehat{B}_j \delta ^T))=c_1 {\rm det}(\delta) {\rm atr}(\widehat{B}_j)={\rm atr}(\widehat{B}_i).\]
Taking into account the inequalities $c_1\neq 0$, ${\rm atr}(\widehat{B}_j) \neq 0$, ${\rm atr}(\widehat{B}_i) \neq 0$, we conclude that ${\rm det}(\delta)\neq0$, i.e. $\delta$ is an invertible matrix

From the equalities $(6)$ and $(7)$ we have \[\gamma \alpha^T+ \delta \widehat{B}_j\beta^T=\gamma \alpha^T +\delta \widehat{B}_j^T \beta^T\]
or \[\delta (\widehat{B}_j-\widehat{B}_j^T)\beta^T=O_{2,n}.\]
From the last equality and from the fact that $\delta$ and $(\widehat{B}_j -\widehat{B}_j^T)$ are invertible (since $\widehat{B}_j$ is a nonsymmetrical matrix of the size $2\times2$), we have $\beta^T=O_{2,n}$.
Then the equality $(5)$ has the form $c_1 I_{n-2}\alpha \alpha^T=I_{n-2}$ and it means that $\alpha$ is an invertible matrix. From $(6)$ we have $\gamma^T=O_{n,2}$. And finally, from equality $(8)$ we have \[\widehat{B}_i=\delta \widehat{B}_j c_1 I_2\delta^T,\]
that runs counter to the choice of the matrices $\widehat{B}_1, \widehat{B}_2, ...$. So, there does not exist a matrix $C$, such that $B_i=C\gamma(C)B_jC^T$ for different $i$ and $j$.

 Since ${\rm det}(B_i)={\rm det}(\widehat{B}_i)$, then $B_i$ is the required set of matrices. It means, that $R(\varphi)=\infty$
\begin{flushright}
$\square$
\end{flushright}

The next result of this paper is

\begin{ttt} Let $G$ be the general linear group ${\rm GL}_n(K)$ or the special linear group ${\rm SL}_n(K)$, where $n\geq 3$ and $K$ is an integral domain. If $K$ contains the ring of integers as a subring and the automorphism group of $K$ is torsion group, then $G$ possesses the $R_{\infty}$-property.
\end{ttt}
\textbf{Proof.}
 Since ${\rm Aut}K$ is a torsion group, then for any ring automorphism $\overline{\delta}$ there is a natural number $m$ such that $ \overline{\delta}^m=id$. Any automorphism $\varphi$ of group $G$ is a standart automorphism, i.e. is a product of inner, central, ring and contragradient automorphisms. As in the Theorem 1 we can rewrite any automorphism $\varphi$ in the form
 \[\varphi = \varphi_D\Lambda^r\Gamma\overline{\delta}.\]
 Also as in Theorem 1, we can consider $\varphi = \Lambda^r\Gamma\overline{\delta},$ and find an infinite set of matrices with equal determinants $det(D)=d$.
 We have two cases that depend on evenness of $r$.
 \par {\it Case} 1. $r$ is an even number. Then $\varphi=\Gamma\overline{\delta}$. Matrices $X$ and $Y$ in $G$ are $\varphi$-conjugated if there is a matrix $Z$ in $G$, which satisfies the equality
 \begin{equation}
 X=ZY\varphi(Z^{-1})=ZY\Gamma\overline{\delta}(Z^{-1})=ZY\gamma(\overline{\delta}(Z^{-1}))\overline{\delta}(Z^{-1}).
\end{equation}
To prove that the number $R(\varphi)$ is infinite, we need to present an infinite set of different $\varphi$-conjugacy classes. According to the equality (9) it is sufficient to find such an infinite set of matrices $A_1, A_2,...$ in ${\rm GL}_n(K)$ with the determinant $d$, such that for any differen $i$ and $j$ there does not exist a matrix $C\in G$, which satisfies the equality
\begin{equation}
A_i=CA_j\gamma(\overline{\delta}(C^{-1}))\overline{\delta}(C^{-1}).
\end{equation}

  Let $X$ be in ${\rm GL}_{n-1}(K)$ and $x$ is an element of the ring $K$. We will use a symbol $X(x)$ to denote a matrix of the form $$X(x)=\begin{pmatrix}
 X &0 \\
 0&x\\
 \end{pmatrix}.$$  Since $K$ contains $\mathbb{Z}$ as a subring, then we find a required infinite set of matrices in the form $A_i(d)$, where $A_i$ belongs to  ${\rm SL}_{n-1}(\mathbb{Z})$. If the following matrices
 $$X(d)=\begin{pmatrix}
 X &0 \\
 0&d\\
 \end{pmatrix},~~~~~~~~~~~~~ Y(d)=\begin{pmatrix}
 Y &0 \\
 0&d\\
 \end{pmatrix},$$
  where $X$ and $Y$ are in ${\rm SL}_{n-1}(\mathbb{Z})$, are $\varphi$-conjugated, then there exists a matrix $Z \in G$ that satisfies the equality
 \begin{equation}X(d)=ZY(d)\gamma(\overline{\delta}(Z^{-1}))\overline{\delta}(Z^{-1}).
\end{equation}
 Since $X$ and $Y$ belong to ${\rm SL}_{n-1}(\mathbb{Z})$ and the group ${\rm Aut}~\mathbb{Z}$ is trivial, we conclude that $\overline{\delta}(X)=X, ~ \overline{\delta}(Y)=Y$ for any ring automorphism $\overline{\delta}$, and it means that $\overline{\delta}(X(d))=X(\delta(d)), ~ \overline{\delta}(Y(d))=Y(\delta(d))$. In consideration of it, operating with the powers of the automorphism $\overline{\delta}$ on the equality (11), we have
\begin{eqnarray}
\nonumber  X(d) &=& Z Y(d)\gamma(\overline{\delta}(Z^{-1}))\overline{\delta}(Z^{-1}), \\
\nonumber  X(\delta(d)) &=& \overline{\delta}(Z)Y(\delta(d))\overline{\delta}(\gamma(\overline{\delta}(Z^{-1})))\overline{\delta}^{2}(Z^{-1}),\\
\nonumber  X(\delta^2(d)) &=& \overline{\delta}^{2}(Z)Y(\delta^2(d))\overline{\delta}^{2}(\gamma(\overline{\delta}(Z^{-1})))\overline{\delta}^{3}(Z^{-1}), \\
\nonumber  \vdots &\vdots& \vdots \\
\nonumber  X(\delta^{m-1}(d)) &=& \overline{\delta}^{m-1}(Z)Y(\delta^{m-1}(d))\overline{\delta}^{m-1}(\gamma(\overline{\delta}(Z^{-1})))Z^{-1}. \end{eqnarray}

If we multiply all these equalities, we have:\[X^m(\widetilde{d})=ZTY^m(\widetilde{d})Z^{-1},\]
where  $$T=\gamma(\overline{\delta}(Z^{-1}))\overline{\delta}(\gamma(\overline{\delta}(Z^{-1})))
\overline{\delta}^{2}(\gamma(\overline{\delta}(Z^{-1})))\dots\overline{\delta}^{m-1}(\gamma(\overline{\delta}(Z^{-1})))=diag(t,\dots,t),$$ $$\widetilde{d}=d\delta(d)\dots \delta^{m-1}(d).$$ As in the first case of Theorem 1, we conclude that the necessary condition of matrices $X(d)$ and $Y(d)$ be $\Gamma\overline{\delta}$-conjugated is the following equality \[(tr(X^m(\widetilde{d})))^n=(tr(Y^m(\widetilde{d})))^n.\]

Now we shall construct a required set of matrices $A_1, A_2, ...$. We examine a matrix $P_x=\begin{pmatrix}
x & 1 \\
-1 & 0 \\
\end{pmatrix}$ and prove that for any natural number $m$ the trace of the matrix $P_x^m$ is a polynomial of the power $m$ (it is important for us that the power of this polynomial is not equal to zero). To do it we prove that there is a polynomial of the power $m$ in the place $(1,1)$ in matrix $P_x^m$, and there are polynomials of smaller powers in the places $(1,2), (2,1), (2,2)$. We use induction on $m$ to prove it. When $m=1$ this statement is obvious. Let $P_x^{m-1}$ has the form $$P_x^{m-1}=\begin{pmatrix} f(x) & g(x) \\
h(x) & p(x) \\
\end{pmatrix}, $$ where $f(x)$ is a polynomial of the power $m-1$ and $g(x),h(x),p(x)$ are polynomials of the powers $r,l,s<m-1$ respectively. Then $$P_x^m=P_x^{m-1}\cdot P_x=\begin{pmatrix}
f(x) & g(x) \\
h(x) & p(x) \\
\end{pmatrix}\cdot\begin{pmatrix}
~x & 1 \\
-1 & 0 \\
\end{pmatrix}=\begin{pmatrix}
xf(x)-g(x) & f(x) \\
xh(x)-p(x) & h(x) \\
\end{pmatrix}.$$
           Now the polynomial in the place (1,1) has a power $m-1+1=m$, in the place (1,2) the power $m-1<m$, in the place (2,1) the power not more than $ max({\rm deg}(xh(x)),{\rm deg}(p(x))=max(l+1,s)<max(m-1+1,m-1)=m$, and finally the polynomial in the place (2,2) of the obtained matrix has the power $l<m-1<m$. It also was required. Hence, we have ${\rm deg}({\rm tr}(P_x^m))=m$. We denote $\psi_m(x)={\rm tr}(P_x^m)$ and  examine now the matrices of the following form: $$A_i=\begin{pmatrix}
P_{a_i}& 0 \\
0 & I_{n-3} \\
\end{pmatrix},$$
 where $\{a_i\}_{i=1}^{\infty}$ is an infinite set of the elements of the ring $K$, such as $(\psi_m(a_i)+n-3+\widetilde{d})^n\neq (\psi_m(a_j)+n-3+\widetilde{d})^n $ for $i \neq j$ (this set exists according to Lemma 1). Then the set $A_1(d), A_2(d), ...$ is the required. Indeed, for any $i\neq j$ the following inequality holds: $$(tr(A_i^m(\widetilde{d})))^n = (\psi_m(a_i)+n-3+\widetilde{d})^n\neq (\psi_m(a_j)+n-3+\widetilde{d})^n=(tr(A_j^m(\widetilde{d})))^n.$$
And it also was required
\par {\it Case} 2. $r$ is an odd number. In this case $\varphi=\Lambda \Gamma\overline{\delta}$. The matrices $X$ and $Y$ are in the same $\varphi$-conjugacy class if for any matrix $Z$ in $G$ we have the following equality
\begin{equation}X=ZY\varphi(Z^{-1})=ZY\Lambda\Gamma\overline{\delta}(Z^{-1})=ZY\Lambda\gamma(\overline{\delta}(Z^{-1}))\overline{\delta}(Z^{-1})=ZY\gamma(\overline{\delta}(Z))\overline{\delta}(Z)^T.
\end{equation}

As in the first case, to prove the $R_{\infty}$-property we find an infinite set of matrices $A_1, A_2, ...$ with the given determinant $d$, where $A_i \nsim_{\varphi} A_j$ if $i \neq j$. We find it in the form $$A_i(d)=\begin{pmatrix}
 A_i &0 \\
 0&d\\
 \end{pmatrix},$$
 where $A_i$ is in ${\rm SL}_{n-1}(\mathbb{Z})$.
We find now the necessary condition for the matrices $X(d)$ and $Y(d)$ of such form to be $\varphi$-conjugated. It means that there exists such a matrix $Z$ that satisfies
\begin{equation}X(d)=ZY(d)\gamma(\overline{\delta}(Z))\overline{\delta}(Z)^T.
\end{equation}

Then, in consideration of $\overline{\delta}(X)=X$, $\overline{\delta}(Y)=Y$, operating with the powers of the automorphism $\Lambda \overline{\delta}$ on the equality (13), we have
\begin{eqnarray}
\nonumber X(d) &=& Z Y(d)\gamma(\overline{\delta}(Z))\overline{\delta}(Z)^T, \\
 \nonumber \Lambda(X)(\delta(d^{-1})) &=&\overline{\delta}((Z^T)^{-1})\Lambda(Y)(\delta(d^{-1}))\overline{\delta}(\gamma(\overline{\delta}(Z^{-1})))\overline{\delta}^{2}(Z^{-1}), \\
 \nonumber X(\delta^2(d)) &=& \overline{\delta}^{2}(Z)Y(\delta^2(d))\overline{\delta}^{2}(\gamma(\overline{\delta}(Z)))\overline{\delta}^{3}(Z)^T, \\
 \nonumber \Lambda(X)(\delta^3(d^{-1})) &=& \overline{\delta}^{3}((Z^T)^{-1})\Lambda(Y)(\delta^3(d^{-1}))\overline{\delta}^3(\gamma(\overline{\delta}(Z^{-1})))\overline{\delta}^{4}(Z^{-1}), \\
 \nonumber \vdots &\vdots& \vdots \\
 \nonumber X(\delta^{2m-2}(d)) &=& \overline{\delta}^{2m-2}(Z)Y(\delta^{2m-2}(d))\overline{\delta}^{2m-2}(\gamma(\overline{\delta}(Z)))\overline{\delta}^{2m-1}(Z)^T, \\
 \nonumber \Lambda(X)(\delta^{2m-1}(d^{-1})) &=&\overline{\delta}^{2m-1}((Z^T)^{-1})\Lambda(Y)(\delta^{2m-1}(d^{-1}))\overline{\delta}^{2m-1}(\gamma(\overline{\delta}(Z^{-1})))Z^{-1}.
\end{eqnarray}
If we multiply all those equalities we have
\[(X\Lambda(X))^m(\widetilde{d})=Z(Y\Lambda(Y))^m(\widetilde{d})TZ^{-1} ,\]
 where $$T=\gamma(\overline{\delta}(Z))\overline{\delta}(\gamma(\overline{\delta}(Z^{-1}))) \dots \overline{\delta}^{2m-2}(\gamma(\overline{\delta}(Z)))\overline{\delta}^{2m-1}(\gamma(\overline{\delta}(Z^{-1})))=diag(t,\dots,t),$$ $$\widetilde{d}=d\delta(d^{-1})\delta^2(d)\delta^3(d^{-1}) \dots \delta^{2m-2}(d)\delta^{2m-1}(d^{-1}).$$ As in the first case it means that \[({\rm tr}((X\Lambda(X))^m(\widetilde{d})))^n=({\rm tr}((Y\Lambda(Y))^m(\widetilde{d})))^n\] is a necessary condition for the matrices $X(d)$ and $Y(d)$ be $\Lambda \Gamma\overline{\delta}$-conjugated .

       We construct now the required set of matrices $A_1, A_2, ...$ that satisfies the condition (12). Let $P_x$ be a matrix of the form $P_x=\begin{pmatrix}
1 & x \\
0 & 1 \\
\end{pmatrix}.$  We show that for any natural number $m$ the trace of the matrix $(P_x\Lambda(P_x))^m$ is a polynomial of the power $2m$. As above, we show that we have a polynomial of the power $2m$ in the place $(1,1)$ in the matrix $(P_x\Lambda(P_x))^m$ and polynomials of the smaller powers in the places $(1,2),(2,1),(2,2)$.  We use a method of induction for it. When $m=1$ we have $$P_x\cdot \Lambda(P_x)=P_x (P_x^{-1})^T=\begin{pmatrix}
1 & x \\
0 & 1 \\
\end{pmatrix} \begin{pmatrix}
~1 & 0 \\
-x & 1 \\
\end{pmatrix}= \begin{pmatrix}
1-x^2 & x \\
-x & 1 \\
\end{pmatrix},$$
i.e. the matrix has a required form.
Let $(P_x\Lambda(P_x))^{m-1}$ have the following form
$$(P_x\Lambda(P_x))^{m-1}=\begin{pmatrix} f(x) & g(x) \\
h(x) & p(x) \\
\end{pmatrix}, $$ where $f(x)$ is a polynomial of the power $2(m-1)$, and $g(x),h(x),p(x)$ are polynomials of the powers $r,l,s<2(m-1)$ respectively. Then

\begin{multline*}(P_x\Lambda(P_x))^m=(P_x\Lambda(P_x))^{m-1} P_x\Lambda(P_x)=\begin{pmatrix}
f(x) & g(x) \\
h(x) & p(x) \\
\end{pmatrix}\begin{pmatrix}
1-x^2 & x \\
-x & 1 \\
\end{pmatrix}=\\
 =\begin{pmatrix}
(1-x^2)f(x)-xg(x) & xf(x)+g(x) \\
(1-x^2)h(x)-xp(x) & xh(x)+p(x) \\
\end{pmatrix}.
\end{multline*}
   We observe now the powers of the obtained polynomials. The polynomial in the place $(1,1)$ has the power $2(m-1)+2=2m$, the polynomial in the place (1,2) has the power $2(m-1)+1=2m-1<2m,$ the polynomial in the place (2,1) has a power which is not greater, than $$max({\rm deg}((1-x^2)h(x)),{\rm deg}(xp(x)))=max(l+2,s+1) <max(2(m-1)+2,2(m-1)+1)=2m,$$ and the polynomial in the place (2,2) has the power $$max({\rm deg}(xh(x)),{\rm deg}(p(x)))=max(l+1,s)<max(2(m-1)+1,2(m-1))=2m-1<2m,$$ as required.  It means that ${\rm deg}({\rm tr}((P_x\Lambda(P_x))^m))=2m$.

  If we denote again $\psi_m(x)={\rm tr}((P_x\Lambda(P_x))^m)$ and examine the following set of matrices
  $$A_i=\begin{pmatrix}P_{a_i} & 0\\
  0& I_{n-3}\end{pmatrix},$$
where $\{a_i\}_{i=1}^{\infty}$ is such set of the elements of the ring $K$, where $(\psi_m(a_i)+n-3+\overline(d))^n \neq \psi_m(a_j)+n-3+\overline(d))^n$ if $i \neq j$ (such set exists according to Lemma 1). Then the set $A_i(d)$ is required. Indeed, for any different $i$ and $j$ the next inequality holds
$$({\rm tr}(A_i(d)\Lambda(A_i(d))))^n=(\psi_m(a_i)+n-3+\overline(d))^n \neq \psi_m(a_j)+n-3+\overline(d))^n=({\rm tr}(A_j(d)\Lambda(A_j(d))))^n,$$
as required. Hence $R(\Lambda\Gamma\overline(\delta))=\infty$.
\begin{flushright}
$\square$
\end{flushright}

\textbf{Remark.} Generally speaking, we could have proved the second case of Theorem 1 using the same way as in the second case of Theorem 2 without using the definition of antitrace. But we wanted to work in terms of matrices $A_i$, without using $A_i\Lambda(A_i)$. Moreover, the application of antitrace gives more elegant proof.

Nasybullov Timur Rinatovich

20/1 Pirogova, Novosibirsk, Russia

phone: 8-923-183-2778

e-mail: timur.nasybullov@mail.ru

\end{document}